\documentstyle[11pt]{article}
\begin{document}

\def\b#1{{\bf #1}}
\def\i#1{{\it #1}}

\title{A particular matrix and its relationship with Fibonacci
numbers}
\author{Mario Catalani
\\Department of Economics\\Via Po 53, 10124 Torino, Italy\\
e-mail mario.catalani@unito.it}
\date{}
\maketitle
\begin{abstract}
\small{Determinants and symmetric functions of the eigenvalues of matrices
characterizing stochastic processes with indepedent increments. Relationships
with Fibonacci numbers are derived.}
\end{abstract}

\section{Introduction}
Let us consider
the $n\times n$ symmetric matrix $\b{A}_n$
\[\b{A}_n=\left [\begin{array}{ccccc}1 & 1 & 1 & \cdots & 1\\
1&2&2&\cdots &2\\
1&2&3&\cdots &3\\
\vdots & \vdots &\vdots &\ddots &\vdots \\
1&2&3&\cdots &n\end{array}\right ],\]
that is
$$a_{ij}=\min(i,\,j),\qquad i,j=1,\,2,\,\cdots ,\,n$$
This matrix is, up to a positive scalar, the covariance matrix of a stochastic
process with increments which possess the same variance and are uncorrelated.

\noindent
For future reference consider also the $(n-k+1)\times (n-k+1)$ matrix
$\b{C}_{n,k}$
\[\b{C}_{n,k}=\left [\begin{array}{ccccc}k & k & k & \cdots & k\\
k&k+1&k+1&\cdots &k+1\\
k&k+1&k+2&\cdots &k+2\\
\vdots & \vdots &\vdots &\ddots &\vdots \\
k&k+1&k+2&\cdots &n\end{array}\right ], \quad \mbox{$k$ integer,
$1<k<n$}.\]

\noindent
For general reference on matrix theory see, for example, \cite{mario},
\cite{gantmacher}, \cite{rao}.

\section{Determinants}
Let us consider the following determinants (dimensions $n\times n$)
$$\Delta_n(i_1,\,i_2,\,\ldots,\,i_n)
=\left |\begin{array}{ccccc}i_1&i_1&i_1&\cdots &i_1\\
i_1&i_1+i_2&i_1+i_2&\cdots &i_1+i_2\\
i_1&i_1+i_2&i_1+i_2+i_3&\cdots &i_1+i_2+i_3\\
\vdots&\vdots &\vdots &\ddots &\vdots\\
i_1&i_1+i_2&i_1+i_2+i_3&\cdots &i_1+i_2+\cdots +i_n\end{array}\right |,$$
and
\begin{eqnarray*}
&&\Theta_n(i_1,\,i_2,\,\ldots,\,i_n,\, i_{n+1})
=\\
&&\qquad =\left |\begin{array}{ccccc}i_1&i_1+i_2&i_1+i_2&\cdots &i_1+i_2\\
i_1&i_1+i_2+i_3&i_1+i_2+i_3&\cdots &i_1+i_2+i_3\\
i_1&i_1+i_2+i_3&i_1+i_2+i_3+i_4&\cdots &i_1+i_2+i_3+i_4\\
\vdots&\vdots &\vdots &\ddots &\vdots\\
i_1&i_1+i_2+i_3&i_1+i_2+i_3+i_4&\cdots &i_1+i_2+\cdots +i_{n+1}
\end{array}\right |,
\end{eqnarray*}
where $i_j, \forall\, j$ are integers $\ge 1$. Let us not that the determinant
of $\b{A}_n$ is $\Delta_n(1,\,1,\,\ldots,\,1)$ and that of the matrix
$\b{C}_{n,k}$ is $\Delta_n(k,\,1,\,\ldots,\,1)$.

\medskip
\noindent
We are going to prove that
\begin{equation}
\Delta_n(i_1,\,i_2,\,\ldots,\,i_n)=\prod_{j=1}^ni_j,
\end{equation}
and
\begin{equation}
\Theta_n(i_1,\,i_2,\,\ldots,\,i_n,\, i_{n+1})=i_1\prod_{j=3}^{n+1}i_j.
\end{equation}
We are going to prove the claim using induction.

\noindent
First of all
\begin{eqnarray*}
\Delta_2(i_1,\, i_2)&=&\left |\begin{array}{cc}i_1&i_1\\
i_1&i_1+i_2\end{array}\right |\\
&=&i_1(i_1+i_2-i_1)\\
&=&i_1i_2,
\end{eqnarray*}
and
\begin{eqnarray*}
\Theta_2(i_1,\, i_2,\, i_3)&=&\left |\begin{array}{cc}i_1&i_1+i_2\\
i_1&i_1+i_2+i_3\end{array}\right |\\
&=&i_1(i_1+i_2+i_3-i_1-i_2)\\
&=&i_1i_3.
\end{eqnarray*}
Now assume that the claim holds for $n-1$. If we evaluate
$\Delta_n(i_1,\,i_2,\,\ldots,\,i_n)$ expanding with respect to the elements
of the first row we have that all the algebraic complements except those
relative to the first two elements of the row are equal to zero: indeed
all the above have identical the first two columns, and we obtain
\begin{eqnarray*}
\Delta_n(i_1,\,i_2,\,\ldots,\,i_n)
&=&i_1
\Delta_{n-1}(i_1+i_2\,i_3,\,\ldots,\,i_n)
-i_1\Theta_{n-1}(i_1,\,i_2,\,\ldots,\,i_n)\\
&=&i_1(i_1+i_2)i_3\cdots i_n-i_1i_1i_3\cdots i_n\\
&=&i_1i_3i_4\cdots i_n(i_1+i_2-i_1)\\
&=&\prod_{j=1}^ni_j.
\end{eqnarray*}
Now reverting to
$\Theta_n(i_1,\,i_2,\,\ldots,\,i_n,\, i_{n+1})$
the same reasoning as before leads to
\begin{eqnarray*}
\Theta_n(i_1,\,i_2,\,\ldots,\,i_n,\, i_{n+1})
&=&i_1
\Delta_{n-1}(i_1+i_2+i_3\,i_4,\,\ldots,\,i_n,\, i_{n+1})\\
&&\qquad -(i_1+i_2)\Theta_{n-1}(i_1,\,i_2+i_3,\,i_4,\,\ldots,\,i_{n+1})\\
&=&i_1(i_1+i_2+i_3)i_4\cdots i_{n+1}-(i_1+i_2)i_1i_4\cdots i_{n+1}\\
&=&i_1i_4\cdots i_{n+1}(i_1+i_2+i_3-i_1-i_2)\\
&=&i_1\prod_{j=3}^{n+1}i_j.
\end{eqnarray*}
This completes the proof.

\noindent
As a corollary we see that immediately we have
$\vert \b{A}_n\vert =1$, and $\vert \b{C}_{n,k}\vert =k$.

\section{Symmetric Functions of Eigenvalues}
Just to fix notation, given $n$ numbers
$$\{a_i\}_1^m=\{a_1,\,a_2,\,\ldots ,\,a_n\}$$
let $S_k^n\left <a_i\right >$ denote the symmetric functions, that is
$$S_k^n\left <a_i\right >=\sum_{\stackrel{i_1,\ldots ,i_k}{i_1<i_2<\cdots <i_k}}
a_{i_1}a_{i_2}\cdots a_{i_k}, \qquad 1\le k\le n.$$
Now given a $n\times n$ matrix \b{A} with eigenvalues $\{\lambda_i\}$ we have
\begin{equation}
\label{eq:sommadeterminanti}
S_k^n\left <\lambda_i\right >=\sum_i\left\vert
\b{A}_{(i)}^{(n-k)}\right \vert,
\end{equation}
where
$$\b{A}_{(i)}^{(h)},\qquad  h=0,\,1,\, \ldots ,\, n-1$$
denotes one of the principal submatrices of \b{A} obtained deleting $h$ rows
and the corresponding $h$ columns. The number of these submatrices is
$${n\choose h}.$$
Note that the submatrices in Equation~\ref{eq:sommadeterminanti} have
dimensions $k\times k$.

\noindent
Considering $\b{A}_n$ we have
\begin{equation}
S_n^n\left <\lambda_i\right > =\vert\b{A}_n\vert = 1,
\end{equation}
\begin{equation}
S_1^n\left <\lambda_i\right > ={\rm tr}(\b{A}_n)={n(n+1)\over 2}.
\end{equation}
A closed inspection of the structure of the determinants in
Equation~\ref{eq:sommadeterminanti} show that
\begin{eqnarray*}
S_k^n\left <\lambda_i\right >&=&\sum_{i_1=1}^{n-k+1}\sum_{i_2=1}^{n-k+2-i_1}
\sum_{i_3=1}^{n-k+3-i_1-i_2}\cdots \\
&&\qquad \cdots \sum_{i_k=1}^{n-i_1-i_2-\cdots -i_{k-1}}
\Delta_k(i_1,i_2\,i_3,\,\ldots,\,i_k)\\
&=&
\sum_{i_1=1}^{n-k+1}\sum_{i_2=1}^{n-k+2-i_1}
\sum_{i_3=1}^{n-k+3-i_1-i_2}\cdots \\
&&\qquad \cdots \sum_{i_k=1}^{n-i_1-i_2-\cdots -i_{k-1}}i_1i_2\cdots i_k.
\end{eqnarray*}
We can develop this formula according the index $i_1$ in the following way
\begin{eqnarray*}
S_k^n\left <\lambda_i\right >&=&1\times \sum_{i_2=1}^{n-k+1}
\sum_{i_3=1}^{n-k+2-i_2}\cdots
\sum_{i_k=1}^{n-1-i_2-\cdots -i_{k-1}}i_2i_3\cdots i_k\\
&&\qquad +2\times \sum_{i_2=1}^{n-k}
\sum_{i_3=1}^{n-k+1-i_2}\cdots
\sum_{i_k=1}^{n-2-i_2-\cdots -i_{k-1}}i_2i_3\cdots i_k\\
&& \qquad +\cdots\,\cdots\\
&&\qquad +(n-k+1)\times \sum_{i_2=1}^1
\sum_{i_3=1}^{2-i_2}\cdots
\sum_{i_k=1}^{k-1-i_2-\cdots -i_{k-1}}i_2i_3\cdots i_k.
\end{eqnarray*}
This shows that we can write
\begin{equation}
S_k^n\left <\lambda_i\right >=\sum_{i=1}^{n-k+1}
iS_{k-1}^{n-i}\left <\lambda_i\right >.
\end{equation}
In the above formula we set $S_0^n\left <\lambda_i\right >=1$.

\noindent
Using this equation we have
$$S_k^n\left <\lambda_i\right >=S_{k-1}^{n-1}\left <\lambda_i\right >
+2S_{k-1}^{n-2}\left <\lambda_i\right >+3S_{k-1}^{n-3}\left <\lambda_i\right >
+\cdots +(n-k+1)S_{k-1}^{k-1}\left <\lambda_i\right >,$$
and
$$
S_k^{n-1}\left <\lambda_i\right >=S_{k-1}^{n-2}\left <\lambda_i\right >
+2S_{k-1}^{n-3}\left <\lambda_i\right >+3S_{k-1}^{n-4}\left <\lambda_i\right >
+\cdots +(n-k)S_{k-1}^{k-1}\left <\lambda_i\right >.$$
It follows
\begin{equation}
\label{eq:ricorsivaesse}
S_k^n\left <\lambda_i\right >=S_k^{n-1}\left <\lambda_i\right >
+\sum_{i=1}^{n-k+1}S_{k-1}^{n-i}\left <\lambda_i\right >.
\end{equation}
This formula makes it easy to evaluate $S_k^n\left <\lambda_i\right >$
for all $n$ and $k$.

\medskip
\noindent
Consider the following identity relating to binomial coefficients
\begin{equation}
\label{eq:ricorsivabinomiale}
{n+k\choose n-k}={n+k-1\choose n-k-1}+\sum_{i=1}^{n-k+1}
{n+k-1-i\choose n-k+1-i}, \qquad 0\le k\le n.
\end{equation}
If we set
$$S_k^n= {n+k\choose n-k},$$ we see that
$$S_k^{n-1}={n+k-1\choose n-k-1},$$
and
$$S_{k-1}^{n-i}={n+k-1-i\choose n-k+1-i}.$$
In this way Equation~\ref{eq:ricorsivaesse} follows from
identity~\ref{eq:ricorsivabinomiale}, and conversely. So we can
conclude that
\begin{equation}
S_k^n\left <\lambda_i\right >={n+k\choose n-k}.
\end{equation}
Note that
$${n+k\choose n-k}={n+k\choose 2k}.$$
Then Identity 1.76 in \cite{gould} says
$$\sum_{k=0}^n{n+k\choose n-k}=F_{2n+1},$$
where $F_i$ is the $i$-th Fibonacci number. It follows
\begin{eqnarray}
\sum_{k=1}^n S_k^n\left <\lambda_i\right >+1&=&
\sum_{k=1}^n {n+k\choose 2k}+1\nonumber\\
&=&\sum_{k=0}^n {n+k\choose 2k}\nonumber\\
&=&F_{2n+1}.
\end{eqnarray}
We can get another recurrence relationship. Indeed from
\begin{eqnarray*}
{n+k\over n-k}{n+k-1\choose n-k-1}&=&
{(n+k)(n+k-1)!\over (n-k)(n-k-1)!(2k)!}\\
&=&{(n+k)!\over (2k)!(n-k)!}\\
&=&{n+k\choose n-k},
\end{eqnarray*}
we get
\begin{equation}
S_k^n\left <\lambda_i\right >={n+k\over n-k}S_k^{n-1}\left <\lambda_i\right >,
\qquad k=1,\,2,\,\ldots ,\,n-1
\end{equation}

\end{document}